\documentclass[12pt,fleqn]{article}  
\usepackage{amsmath,amssymb,amsfonts}  

\newcommand{\R}{\mathbb{R}} 
\newcommand{\C}{\mathbb{C}}

\begin{document}   


\voffset=-.5 truecm 
\textheight=8.25 truein

\title{A new integrable system on the sphere}   
\date{}   
   
\author{%
Holger R. Dullin%
\thanks{Department of Mathematical Sciences, Loughborough University, LE11 3TU  UK;
h.r.dullin@lboro.ac.uk}, 
Vladimir S. Matveev%
\thanks{Mathematisches Institut, Universit\"at Freiburg, 79104 Germany
 matveev@email.mathematik.uni-freiburg.de}}   
\maketitle   
   
\begin{abstract}   
We present a new Liouville-integrable natural  Hamiltonian system   
on the (cotangent bundle of the) sphere $S^2$. 
The second integral is cubic in the momenta.   
\end{abstract}   
   
{MSC2000: 37J35, 58F07, 58F17, 70H06, 70E40}

\section{Introduction}   

A Hamiltonian system is called natural if 
its Hamiltonian is a sum of a positive-definite 
kinetic energy and a potential.  
Natural Hamiltonian systems on (cotangent bundles of) closed  surfaces  
admitting  integrals polynomial in momenta are  
interesting for a number of reasons:  
  
\begin{enumerate}  
\item They are classical:  The formulation  
of the problem is due at least to Darboux \cite{Darboux}.

\item If a natural Hamiltonian system admits a real-analytic integral,  
then it admits an integral polynomial in momenta \cite{Whittaker}.  
  
\item In particular, all known natural integrable Hamiltonian  
systems on surfaces have integrals polynomial in momenta.  
  
\item The existence of the integral polynomial in momenta of degree  
one or  two  has   a very clear geometric background: the  
existence of an integral of degree one implies the existence of a  
one-parametric family of symmetries of the system.  The existence  
of an integral of degree two implies the existence of so-called  
separating variables.  
\end{enumerate}

Natural Hamiltonian systems on closed surfaces that admit a  
nontrivial integral which is polynomial in momenta of degree one  
or two are completely understood: there exists  a complete  
description and classification, see, for example, \cite{BF}.  
Up to now, only one family of natural  
Hamiltonian systems on a closed surface admitting an integral of  
degree three and admitting no (nontrivial) integral of degree one  
or two is explicitly known, namely the family of  system obtained 
from the  Goryachev-Chaplygin top.  
We now present a second, distinct, family of natural Hamiltonian 
systems with a cubic integral.

\subsection{The result} 
Consider the sphere $S^2 \subset \R^3$ of radius 1 with  
the spherical coordinates $(x,y,z) = (-\sin\theta\cos\phi,-\sin\theta\sin\phi,
 \cos\theta)$ and the following two functions  
$H$ (the Hamiltonian) and $F$ (the second integral) 
on the cotangent bundle of the sphere without poles $z={\pm}1$.  
Let $A, c, s \in \R$ be parameters with $s > 1$ and define  
\[
     W(z)  =  z + s, \quad
      P(z) = 3z^2 + 4 s z + 1, \quad
      Q(z) = 3z^2 + 2 s z - 1
\]
and 
\[ 
    G(z)  =  \frac{P(z)}{(2 W(z))^2}. 
\]
Then, the function $H$ is given by  $H :=  K+ V$, where  
\begin{align*}
K & :=\frac12 \left( \frac{1}{\sin^2\theta} +{G(\cos\theta)}\right) p_\phi^2 + \frac12 p_\theta^2 \\
V & := A\frac{ \sin\theta}{ \sqrt{W(\cos\theta)} } \cos \phi+ \frac{c}{W(\cos\theta)},
\end{align*}
and $F$ is defined as
 \[  
 F := 2Hp_\phi - p_\phi^3 +
A  \cos(\phi)\frac{Q(\cos\theta)}{\sqrt{W(\cos\theta)}\sin\theta}
p_\phi + 2A\sin \phi \sqrt{W(\cos\theta)}  p_\theta   \,.
\]

\vspace{2ex}
 \noindent {\bf Proposition: }  
 {\em  The functions $H$ and $F$ can be analytically continued to  
the cotangent bundle of the whole sphere.  The continuation is also  
polynomial in momenta (of degree 2 for H and of degree 3 for $F$).}

\vspace{2ex}  We will denote the continuations of $H$, $K$, $F$  
  by the same letters $H$, $K$, $F$; in particular in the  
theorem  below we mean the continued functions
defined on the cotangent bundle of the whole sphere.  
  
\vspace{2ex}  
  
\noindent {\bf Theorem:} 
{\it The following statements hold:  
\begin{enumerate}
  
\item   The functions $H$ and  $F$ commute with respect  to the standard Poisson bracket on $T^*S^2$
and  are functionally  independent  
  
\item  The  kinetic energy $K$ is positive definite  
  
\item  If $A\ne 0$, the Hamiltonian $H$ does not admit a (smooth) nontrivial   integral which is  polynomial in velocities of degree less than three and which is linearly independent of  
$H$.   
  
\end{enumerate}  
}  
  
\vspace{2ex}  
  
In other words,  we found a new 1) Liouville integrable,
2) natural Hamiltonian system on the sphere with 3) an integral cubic in momenta.

\subsection{Why the system is new} \label{sec:new}  
  
The previously known family comes  from  the Goryachev-Chaplygin top 
using symplectic reduction, see \cite{BKF} for the details. 
It is a natural Hamiltonian system on the sphere. In the standard  
spherical  coordinates its Hamiltonian $H_1$ is $K_1+V_1$, where  
the kinetic energy $K_1$ and the potential $V_1$ are 
\begin{align*}
K_1 & :=\frac{1}{2}\left(\frac{\cos^2\theta}{\sin^2\theta}+4
\right)p_\phi^2+ \frac{1}{2}p_\theta^2  \\
V_1 & := A_1\sin \theta \sin \phi\, ,
\end{align*}
see, for example, \cite{DMT}.  
The integral $F_1$ for the Hamiltonian $H_1$ is given by  
$$
F_1:= H_1p_\phi-2p_\phi^3-\frac{A_1}{2}
p_\theta\cos(\phi)\cos(\theta) +\frac{A_1}{2}{\frac {{}\,\sin \phi
\left( 3\, \left( \cos  \theta  \right) ^{2}-2 \right) }{\sin
\theta }}p_\phi. $$

It is easy to see  that the system we found is essentially
different from the Goryachev-Chaplygin system, i.e.\ that there
exists  no diffemorphism $D:S^2\to S^2$ such that $D_*\, K= \alpha
\, K_1$ and $D_*\,  V =\beta V_1$ for certain constants $\alpha$,
$\beta$.

Consider the metrics $g,$ $g_1$ corresponding to the kinetic  
energies $K$, $K_1$, respectively. They  are given by the  
formulae:  
\begin{eqnarray} \label{eqn:metric}
ds^2 &=&  d\theta^2 + \frac{d\phi^2}{1/\sin^2\theta + G(\cos
\theta)} \label{metric1}
\\ \label{eqn:metric1}
ds_1^2&=&{d\theta^2}+ \frac{d\phi^2}{4+\cot^2\theta}
   = {d\theta^2}+ \frac{d\phi^2}{1/\sin^2\theta + 3}
\label{metric2}.
\end{eqnarray}
  
By direct calculation,  it is easy to see that level lines of the  
curvature of both metric are the lines $\{\theta=\textrm{const}\}$.  
More precisely, since the formulae above do not depend on $\phi$,  
the curvature is constant along every such line. It is easy to see  
that the curvature is not constant.  
  
Then, the diffeomorphism $D$, if it exists, must take the lines  
$\{\theta=\textrm{const}\}$ of the first system to the lines  
$\{\theta=\textrm{const}_1\}$ of the second system. Since the  
lines $\{\theta=\textrm{const}\}$ are orthogonal to the lines  
$\{\phi=\textrm{const}\}$ in the first and in the second metric,  
the diffeomorphism $D$  must also take the lines  
$\{\phi=\textrm{const}\}$ of the first system to the lines  
$\{\phi=\textrm{const}_1\}$ of the second system. Using that  the  
coordinate $\phi$ is periodic with the period $2\pi$, we obtain  
that  the diffeomorphism $D$ must be given by formulae $$  
\theta_{new}=\theta_{new}(\theta), \ \ \ \phi_{new}=\pm \phi+  
\phi_0 \, .$$  
 Then, the pull-back of the metric $g_1$ has the  
form $$ \left(\frac{\partial \theta_{new}}{\partial  
\theta}\right)^2 d\theta^2  + \frac{d\phi^2}{4 + \cot^2  
\theta_{new}}.$$ Then, $\frac{\partial \theta_{new}}{\partial  
\theta}=\pm 1$, and the  diffeomorphism $D$ must be given by the  
formula  
 $$  
\theta_{new}=\theta \ \ \ \ \ \mathrm{or} \ \  \theta_{new}=\pi -  
\theta, \ \ \ \ \ \phi_{new}=\pm \phi+ \phi_0 \,  . $$ Clearly,  
the pull-back of the metric $g_1$  with the help of such  
diffeomorphism  is not proportional to $g$. Otherwise $G$ would
need to become a constant. But then $P$ would have a double root at $z = -s$,
which is impossible.
Thus, our system is essentially different from the reduced Goryachev-Chaplygin top.

Selivanova \cite{Selivanova} proved the existence of an additional
family of natural Hamiltonian systems admitting  integrals of degree 3 in  
momenta.  This family is not explicit. Instead, a nonlinear  
differential equation is derived and it is proved that certain  
solutions of this equation allow to construct a geodesic flow  
admitting an integral of degree three and admitting no integral of  
smaller degree.

At present we do not know whether our system overlaps  
with the systems considered by Selivanova.  Our guess is that this 
is not the case.  
Here are some reasons supporting this claim:  
\begin{enumerate}  
\item It is not clear whether the  
system of  Selivanova is real-analytic; she claims the smoothness only.  
  
\item Selivanova's family  contains (the reduction of) Goryachev-Chaplygin system; 
our family does not.  
  
\item The ansatz of Selivanova is more restrictive than the form of our system.  
  
\end{enumerate}

\section{Proof of the proposition}
  
\subsection{Analytic proof} 
\label{proposition}  
Clearly, the functions $V,K,F$  depend real-analytically  
on $(\phi,\theta)$ and on the corresponding momenta $p_\phi, p_\theta$,  
so the only possible problem could appear near the points $z=\pm  
1$. We consider $(x,y)$  as a local coordinate system near these  
points  of the sphere. The coordinates $(x,y)$ are related to  
coordinates $(\phi,\theta)$ by the following real-analytic change of variables:  
\begin{align*}  
   \sin\theta &= \sqrt{x^2+y^2}, &  p_\theta &= (x p_x + y p_y) \sqrt{(x^2+y^2)^{-1}-1} \\  
   \tan\phi &= y/x, \quad             & p_\phi &= x p_y - y p_x  \,.  
\end{align*}  
In the new local coordinates the kinetic energy and potential are  
\[  
   K = \frac12\left( p_x^2 + p_y^2 - (xp_x + yp_y)^2 + (xp_y-yp_x)^2 G(z) \right) ,  
\]  
\[  
    V(x,y) = -A \frac{x}{\sqrt{W(z)}} + \frac{c}{W(z)}, \ \ \ \textrm{where} \quad z = \pm \sqrt{1-x^2-y^2}  
    \,.  
\]  
The sign of $z$ depends on whether a neighbourhood of the north or  
south pole is considered. Clearly, $K$ and $V$ are real-analytic.  
Thus, the Hamiltonian   $H$ is a real-analytic  function on the  
whole $T^*S^2$. 
Since we will need it for the proof of  the second statement of  the Theorem, 
we remark that the kinetic energy is positive definite near the poles, 
where $x$ and $y$ are small.

Now let us show that the second integral also is a real-analytic function  
on the whole  $T^*S^2$. The first two terms of $F$  are  
real-analytic,  since the momentum $p_\phi$ is $x p_y - y p_x$ in  
the new coordinate system, and therefore is real-analytic, and we  
already have proven that $H$ is real-analytic. The sum of the  
third  and the fourth terms in the new coordinate system is 
$$ 
A\left(\frac{-xy}{\sqrt{\sqrt{1-x^2-y^2}+s}}p_x +  
\frac{2-3x^2-2y^2+2s\sqrt{1-x^2-y^2}}{\sqrt{\sqrt{1-x^2-y^2}+s}}p_y  
\right), $$ 
and this is clearly real-analytic.

\subsection{Geometric proof}
  
The second integral $F$ is cubic in the momenta, like the  
Goryachev-Chaplygin top. This top does not appear as a limiting  
case in our family, but is serves as a motivation for the  
following coordinate transformation, which uses global coordinates  
$(x,y,z)$ on a sphere embedded in $\R^3$. The angular momentum  
$(L_x,L_y,L_z)$ where $L_z = x p_y - y p_x = p_\phi$ etc.\ with  
cyclic permutations, satisfies the usual Poisson structure on the  
sphere with non-vanishing brackets  
\[  
  \{ L_x, L_y \} = L_z, \quad 
   \{x, L_y \} = z, \, \{ z, L_y \} = -x  
\]  
and cyclic permutations thereof. This bracket has  
Casimirs $x^2+y^2+z^2=1$ and $xL_x + yL_y + zL_z=0$.  
The global Hamiltonian is $H = K+V$, where  
\[  
   K  = \frac12 \left( L_x^2 + L_y^2 + (1+G(z))L_z^2 \right)  
\]  
\[  
      V(x,y,z) = -A \frac{ x}{\sqrt{W(z)}} + \frac{c}{W(z)},  
\]  
and is clearly real-analytic. The integral reads  
\[  
\begin{aligned}  
  F & =  
2  H L_z -  L_z^3 -  
A \frac{Q}{\sqrt{W}}\frac{xL_z}{1-z^2}  
   - 2A\sqrt{W} \frac{y (xL_y-yL_x)}{1-z^2}  \\
   & = 2  H L_z -  L_z^3   
+ \frac{A}{\sqrt{W}}( x L_z + 2W L_x )
   \,.  
\end{aligned}  
\]
In the second line the apparent singularity in $F$ is removed by using
both Casimirs $x^2+y^2+z^2 = 1$ and $x L_x + y L_y + z L_z = 0$ and 
in addition the identity $Q = 2zW + z^2 - 1$.  
The integral is cubic in the momenta, and clearly is real-analytic. 

Like the Goryachev-Chaplygin  
top the system is only integrable on the level set of the Casimir  
$x L_x + y L_y + z L_z = 0$. Writing the Hamiltonian in this form  
suggests an interpretation of a spinning top with zero  
angular momentum (the vanishing Casimir), an orientation
dependent moment of inertia, and a non-standard potential. 
Note, however,  
that $1+G$ is not positive for all $z \in [-1,1]$ for certain  
values of $s$, which means that a moment of inertia may 
pass through infinity for such value of $s$.
Nevertheless, the kinetic energy in $H$ is positive  
definite for all $s$, see Section~\ref{3.1}. 
This is  
possible because even when $1+G < 0$ because $x L_x + y L_y + z  
L_z = 0$.

Introducing $\xi = x + i y$ and $\eta = L_x + i L_y$ and writing $V(x,y,z) = x U(z)$ 
the equations of motion are
\begin{align}
  -i \dot \xi    & = - L_z  ( 1 + G) \xi + z \eta  \\
  -i \dot \eta & =  - \frac12 G' L_z^2 \xi - L_z G \eta + z U - \frac12 U' \,(\xi^2 + \xi \bar \xi) \,.
\end{align}  
Here $z^2 = 1 - \xi \bar \xi $ and $-2 z L_z = \xi \bar \eta + \eta \bar \xi$,  
or the equations are completed with 
$ \dot z  = x L_y - y L_x$ and $ \dot L_z  =  -y U$.

\section{Proof of the Theorem}  
  
\subsection{Liouville integrable}

It is tedious but easy to check that the functions $H$ and $F$ commute, 
since they are given by explicit formulae. 
The calculation goes as follows.
The canonical bracket is a quadratic polynomial in the momenta. The term independent 
of the momenta vanishes because $Q = 2z W - 1 + z^2$.
The coefficient of $p_\theta p_\phi$ vanishes because $P = 4 z W + 1 - z^2$. 
The nontrivial fact is that the coefficients of $p_\phi^2$ which contains $P$, $Q$ 
and derivatives thereof vanishes.

$H$ and $F$ are independent because $H$ has terms quadratic in momenta 
(and not a square of a linear function of momenta) 
and $F$ has terms that are cubic in momenta
(and not a square of a linear function of momenta).
So $H$ and $F$ have vanishing Poisson bracket and are functionally 
independent, 
hence the system is Liouville integrable.

\subsection{Natural} \label{3.1}  
  
As we already have remarked in Section~\ref{proposition}, the  
kinetic energy $K$ is positive near the poles of the sphere. Thus,  
it is sufficient to show that the metric (\ref{metric1}) is  
positive definite at regular points of the spherical coordinate  
systems, i.e. we have to show that $C:=1/\sin^2\theta + G(\cos  
\theta)$ is positive for $0<\theta<\pi$. Substituting $z$ instead  
of $\cos \theta$, we obtain  
\[  
 C(z)=\frac{1}{4} \frac{6z^{2}+12s z+4 {s}^{2}-3{z}^{4}-4s{z}^{3}-1}{(1-{z}^{2})  (z+s)^{2}} \, . 
\]  
 Since the denominator  
$\left( 1-{z}^{2} \right)  \left( z+s \right) ^{2}$ is always  
positive for $-1<z=\cos(\theta)<1$, we need to prove that the numerator 
$$  
R(z):={\frac{3}{2}{z}^{2}+3s z+{s}^{2}-\frac{3}{4}{z}^{4}-s{z}^{3}-\frac{1}{4}} 
$$ 
is positive.  
Its value at $z=-1$ is $(s-1)^2$ and hence positive; its derivative is 
$$
\frac{\partial R(z)}{\partial z} = 3(1-z^2)(z+s) \,.
$$ 
Since the derivative is also positive for $-1<z=\cos(\theta)<1$, $R(z)$ is positive  
for $-1<z=\cos(\theta)<1$. 
Hence the metric (\ref{metric1}) is positive definite.

\subsection{Cubic}  
  
 Let us prove that,  for $A\ne 0$, our system does not admit  
 an integral linear in momenta, and that every integral quadratic in momenta is  
 proportional to the Hamiltonian.

 Suppose the function $P_2+P_1+P_0$, where every $P_i$ is a  
 homogeneous polynomial in momenta of degree $i$,  is an integral for $H=K+V$.  
 Then, the Poisson bracket $\{K+V, P_2+P_1+P_0\}$ must vanish. We  have  
 \begin{align} \label{bracket}  
 0 & =  \{K+V, P_2+P_1+P_0\}   \nonumber \\
 & =  \{K,P_2\}+ \{K,P_1\}+   \{K,P_0\} + \{V,P_2\}  +  \{V,P_1\}+  \{V,P_0\} \, .   
 \end{align}  
  
Since the Poisson bracket  of two functions $F_1,F_2:T^*M^2\to  
\mathbb{R}$ is given by $$ \{F_1,F_2\}=\frac{\partial  
F_1}{\partial x}\frac{\partial F_2}{\partial p_x}+ \frac{\partial  
F_1}{\partial y}\frac{\partial F_2}{\partial p_y} -\frac{\partial  
F_1}{\partial p_x}\frac{\partial F_2}{\partial x}-\frac{\partial  
F_1}{\partial p_y}\frac{\partial F_2}{\partial y}, $$ where $x,y$  
are local coordinates on $M^2$ and $p_x,p_y$ are the corresponding  
momenta, every term on the right-hand side of (\ref{bracket}) is  
polynomial and homogeneous in momenta: $\{K,P_2\}$ of degree $3$;  
$\{K,P_1\}$ of degree $2$;  
$\{K,P_0\} +  \{V,P_2\}$ of degree $1$ and $\{V,P_1\}$ of degree $0$. 
The Poisson bracket $\{V,P_0\}$ evidently vanishes.  
 For (\ref{bracket}) to vanish all the other degree 0, 1, 2, and 3 polynomials 
 must also vanish identically.
 
  Let us  prove that the linear term $P_1$ must be zero.  
 Vanishing of $\{K,P_1\}$  implies that $P_1$ is a linear integral  
 of the geodesic flow of the metric given by $K$. In the  spherical  
 coordinates $\phi, \theta$ on the sphere the metric is given by (\ref{eqn:metric}).  
We have already remarked in Section \ref{sec:new} that the 
level curves of the curvature are the curves $\{\theta=\mathrm{const}\}$.  
  
Suppose the integral $P_1$ has the form  
$v_1(\theta,\phi)p_\theta+v_2(\theta,\phi)p_\phi$. Then, by  
Noether's Theorem, the vector field $v:=(v_1,v_2)$ is a Killing  
vector field, i.e. its flow preserves the metric. Then, the flow  
of $v$ preserves the  curvature, which implies that the component  
$v_1$ must be zero.  
  
Now, using that $\{V,P_1\}$ vanishes, we have that $\{V,  
v_2p_\phi\}$ is zero, which in the case $v_2\not\equiv 0$ implies  
$\frac{\partial V}{\partial \phi}=0$, which is definitely not the  
case. Thus, $v_2\equiv 0$,  and $P_1\equiv 0$.

Now let us show  that $P_2+P_0$ is proportional to the Hamiltonian $H$.  
Since $\{K,P_2\}=0$, $P_2$ is an  integral  of the geodesic flow  
for the metric (\ref{eqn:metric}). As we already mentioned, for 
geodesic flows on closed surfaces the  
integrals quadratic in velocities are completely understood. 
In particular, if the  
curvature of the surface is not constant, then  the space of  
quadratic integrals is at most two-dimensional  
\cite{dissertation,Kio}.  
  
The curvature of the metric (\ref{eqn:metric}) is not  constant, and  
the following two linearly independent functions are integrals of  
its geodesic flow: the first is $K$ itself, and the second is  
$p_\phi^2$. Then, for certain constants $\beta,\alpha$, we  
have $P_2=\alpha p_\phi^2 + \beta K$. Our goal is to show that  
$P_2+P_0$ is proportional to  $H$. 
>From $\{K,P_0\} + \{V,P_2\}=0$ we find
$\{K,P_0 - \beta V \} + \alpha  \{V,p_\phi^2\}=0$. 
Hence,  
$$  
 p_\theta \frac{\partial (P_0- \beta V)}{\partial  
\theta}+ p_\phi\frac{1}{\sin^2\theta+G(\cos \theta)}\frac{\partial  
(P_0 - \beta V)}{\partial \phi} -2\alpha p_\phi \frac{\partial V}{\partial  
\phi}=0 \,.  
 $$  
Therefore $\frac{\partial (P_0-\beta V)}{\partial \theta}=0$, and then   
$\frac{\partial (P_0-\beta V)}{\partial \phi}=2\alpha (\sin^2\theta+G(\cos  
\theta))\frac{\partial V}{\partial \phi}$. 
By cross differentiation we find  
$$ 
2\alpha \frac{\partial \left( (\sin^2\theta+G(\cos  
\theta))\frac{\partial V}{\partial \phi}\right)}{\partial  
\theta}=0, 
$$  
which is possible only if $\alpha=0$ or $A=0$. 
If $\alpha = 0$ then $P_0 - \beta V = {\rm const}$, 
and so after absorbing the constant in the Hamiltonian 
we have $\beta = 1$ and the integral $P_2+P_1+P_0$ is proportional to $H$. 
This proves the Theorem.

\section{Geodesic flows with cubic integrals}  
  
A geodesic flow is  a  natural Hamiltonian system whose potential  
energy is identically zero.  
  
According to Maupertuis's principle, an integrable natural  
Hamiltonian system  immediately gives a family of  integrable  
geodesic flow, see \cite{BKF} for details. If the integral of the  
system is polynomial in momenta, the integrals of the geodesic  
flows are  also polynomial of the same degree.

In our case, the Hamiltonian of the geodesic flow is  
 given by $$H_{geod} =  \frac{p_\theta^2 +  
{p_\phi^2}{\left(1/\sin^2\theta + G(\cos  
\theta)\right)}}{2(h-V)}.$$ Here $h$ is a constant such that  
$h$ is larger than the maximum of the potential $V(x)$ on the sphere $S^2$.
The integral of the third degree is  
\begin{equation}\label{integral12}  
\begin{aligned}  
   F_{geod} = & 2Hp_\phi - 2Vp_\phi+ p_\phi^3  
    +H_{geod} \Big( 2Vp_\phi \\ 
  &  +    
 A  \cos(\phi)\frac{Q(\cos\theta)}{\sqrt{W(\cos\theta)}\sin\theta}  
p_\phi + 2A\sin \phi \sqrt{W(\cos\theta)}  p_\theta \Big)  \,.  
\end{aligned}  
\end{equation}  
  
The metric corresponding to $H_{geod} $ is  
 \begin{equation}  
\label{metric11} ds^2_{geod} =\left(h-V\right) \left(d\theta^2 +  
\frac{d\phi^2}{1/\sin^2\theta + G(\cos  
\theta)}\right).
\end{equation}

The investigation of geodesic flows that admit an integral polynomial in momenta 
is a very classical subject. If the degree of integrals is one or two, they  
were  well-understood already in 19th century, see for example  
\cite{Darboux}. There are a lot of local examples of geodesic  
flows admitting integral polynomial in velocities of degree 3, see  
the survey \cite{Hijetarinta} for details, and only very few global, 
i.e.\ on closed or complete surfaces.

Let us recall the known results about geodesic flows on closed  
surfaces admitting an integral polynomial in momenta.  For a more  
detailed review see \cite{BF}.  
  
First of all, in view of results of Kolokoltsov \cite{Kol}, 
a geodesic flow of a surfaces of genus greater than two cannot  
admit a nontrivial integral polynomial in momenta. Then, an  
orientable surfaces admitting such geodesic flows must be the  
sphere or the torus.  
  
We  collect  the main results about existence  
 and classification of  metrics on the torus and  on  
 the sphere  whose geodesic flow  
 admits a nontrivial integral polynomial in momenta in the following table:  
  
\vspace{3ex}  
  
\begin{tabular} {|c|c|c|}  
\hline  
         & Sphere $S^2$  &Torus $T^2$ \\  
         \hline  
Degree 1 & All is known & All is known \\ \hline  
 Degree 2 & All is  
known & All is known\\ \hline  
 Degree 3 & Series of examples & Partial negative results \\  
\hline Degree 4 & Series of examples & Partial negative results \\  
\hline Degree $\ge$ 5 & Nothing is known  & Nothing is known \\  
\hline  
\end{tabular}  
  
\vspace{3ex}  
  
The following notation is used in the table. 
``Degree" means the smallest degree of an integral polynomial in velocities.  
``All is known" means that there exists an effective description and classification 
(can be found in \cite{BMF}).  
  
There exists only one explicit ``Series of examples" for degree  
three. It  comes from the Goryachev-Chaplygin case of rigid body  
motion (by applying symplectic reduction and Maupertuis's  
principle), see \cite{BKF}. The most valuable "partial negative  
results" are due to Byaly\u i \cite{Byalyi} and Denisova and  
Kozlov \cite{Denisova}.  They proved that if a natural Hamiltonian  
system on the torus  whose kinetic energy is given by a flat  
metric admits an integral polynomial in momenta of degree three,  
it admits an integral  linear in momenta.  
  
Thus, our system gives one more series of examples of integrable  
geodesic flows on the sphere whose integral is polynomial in  
momenta of degree 3. Below, in Sections~\ref{4.1}, \ref{4.2}, we  
show that these geodesic flows do not admit an  integral which is  
polynomial in momenta of degree less than 3, and that these  
examples are different from examples coming from  
Goryachev-Chaplygin.   Note, that this does not follow directly from 
the Theorem.

Kiyohara \cite{Kyoh} proved the existence of one more series of  
examples of integrable geodesic flows on the sphere whose integral  
is polynomial in momenta of degree 3. Although his examples are  
not explicit (similar to Selivanova \cite{Selivanova} he writes a  
nonlinear differential equation and claims that its solutions give  
integrable geodesic flows), it is easy to see that his examples  
are different from ours, since they have regions of constant  
curvature and therefore can not be real-analytic.

\subsection{Cubic}
\label{4.1}  
  
We need to prove that the Hamiltonian $H_{geod}$ has no integral which is 
polynomial in velocities of degree 1,2. 
We will prove it assuming that the energy level $h$  is sufficiently large.  
Our proof for arbitrary $h$ is based on the technique  
developed in \cite{BF} and is very long. It will be published  
elsewhere.  
  
It is easy to see that the integral $F_{geod}$ is  not the third  
power of a function linear in momenta and is not the product of a  
function linear in momenta and the Hamiltonian. Then, if  
$F_{additional}\not\equiv 0$  is linear or quadratic  in momenta  
integral for $H_{geod}$, it must be functionally independent of  
$F_{geod}$. Then, the system is resonant. Then, all geodesics are  
closed. The theory of metrics all whose geodesics are closed is  
well-developed, see \cite{besse}. It is known \cite{gromoll} that, on the  
two-sphere, geodesics of such metrics have no self-intersections.  
But if $h$ is very big, certain  geodesics of our metric must have  
self-intersections: indeed, for huge $h$, the geodesics of  our  
metric are very close to the geodesic of the metric  
(\ref{metric1}), and by direct analysis it is easy to check that  
certain  geodesics of (\ref{metric1}) do have self-intersections.
  
Thus, for very big $h$,  the Hamiltonian $H_{geod}$ admits no  
integral of  degree one or two.  
  
\subsection{Why  $ds^2_{geod}$  is new}
\label{4.2}  
  
The examples coming   from Goryachev-Chaplygin are the metrics
\begin{equation}  
\label{metric12} ds^2_{2} =\left(h_1-V_1\right)\left({d\theta^2}+  
\frac{d\phi^2}{4+\cot^2(\theta)}\right)\, .\end{equation} Their  
Hamiltonians are $$ H_2:=\frac{{p_\theta^2}+  
\left({4+\cot^2(\theta)}\right){p_\phi^2}}{2(h_1-V_1)} \, ,
$$  
and the integral for the system with the Hamiltonian $H_2$ is  
$$
\begin{aligned}  
 F_2:=&H_1p_\phi- 2p_\phi^3- V_1p_\phi \\
  &+ \left(V_1p_\phi-\frac{A_1}{2} p_\theta\cos(\phi)\cos(\theta)  
+\frac{A_1}{2}{\frac {{}\,\sin \phi \left( 3\, \left( \cos  \theta  
\right) ^{2}-2 \right) }{\sin \theta }}p_\phi\right) H_{2} \, .
\end{aligned}
$$  
  
Our goal is to show that there is no diffeomorphism $D:S^2\to S^2$  
that takes the metric (\ref{metric11}) to the metric (\ref{metric12})
or a constant multiple thereof.  
Since the potential energy $V_1$  depends linearly  on the parameter $A_1$,
it is sufficient to show that  there exists no diffeomorphism $D:S^2\to S^2$  
that takes the metric (\ref{metric11}) to the metric (\ref{metric12}). 
  
First let us show that such diffeomorphism $D$ must be given by  
the formulae $\theta_{new}=\theta_{new}(\theta)$,  
$\phi_{new}=\pm\phi+\phi_0$, where $\phi_0$ is a constant.  
 It is known, that the Goryachev-Chaplygin system is nonresonant on every energy level, 
 see for example \cite{Orel} where the rotation function is calculated explicitely.  
 Then, the geodesic flow of  (\ref{metric12}) is nonresonant as well. Then,  
 the diffeomorphism $D$ must take the integral (\ref{integral12}) 
 to the integral (\ref{metric12}) or a constant multiple thereof.

It is known that a Riemannian metric on a (oriented) surface  
defines a complex structure on it. That is, in a neighbourhood of  
every point we can find a  coordinate $z=x+iy$ such that the  
metric looks like $\lambda(x,y)(dx^2+dy^2)$. In such  complex  
coordinates,  
an  integral of degree three  
has the form 
\begin{equation}\label{integral} 
a(z)p^3+b(z)p^2\bar p+ \bar b(z)p\bar p^2+\bar a(z) \bar p^3.
\end{equation} 
Here $a,b$ are functions;  
``bar" means  the complex  conjugation, $p$ is the complex momentum
 corresponding to  $z$. The  
functions $a,b$ are  a-priory not assumed to be holomorphic.  
Note that this form is not restrictive:  If  (\ref{integral}) 
is a real-valued  
function, the coefficients at $\bar p^3$ and $p^2\bar p$ must be  
complex-conjugate to the  coefficients at $p^3$ and $\bar p^2 p$,  
respectively.  
  
We will use the following result of Kolokoltsov from \cite{Kol};  
in a weaker form, this result has been already known to Birkhoff  
\cite{Birkhoff}.

\vspace{2ex}  
  
  \noindent {\bf Lemma: }{\it  Suppose (\ref{integral}) is an integral for  
the geodesic flow of $g$ and that it is not the product of the Hamiltonian  
and an integral linear in velocities. Then, \begin{equation}  
\label{Bir} \frac{1}{a(z)}dz^3  
\end{equation}  
is a meromorphic $(3,0)$-form. }  
  
\vspace{2ex}  
  
Evidently, the form (\ref{Bir}) has no zeros. It could have poles.  
In fact, because of Abel's Lemma, it always has poles on the  
sphere.

The metrics (\ref{metric11},\ref{metric12}) and the integrals  
$F_{geod},F_2$ are given explicitly.  Direct calculation of the  
form (\ref{Bir}) shows that the form (\ref{Bir}) for the metric  
(\ref{metric11}) and for the integral  $F_{geod}$ is very similar  
to the form (\ref{Bir}) for the metric (\ref{metric12}) and for the  
integral  $F_{2}$: they have precisely two poles. The poles are  
located at $\theta=0,\pi$.  
Then, the diffeomorpism $D$ takes the points $\theta=0,\pi$ to the points $\theta=0,\pi$.  
  
It is known  that a holomorphic diffeomorphism of the Riemann  
sphere $\bar {\C}$  that preserves the points $z=0$ and $z=\infty$  
is a linear transformation $z\mapsto \alpha z$. A holomorphic  
diffeomorphism that interchange the points is a transformation  
$z\mapsto \frac{\alpha}{ z}$. In both cases, the diffeomorphism  
$D$ has the form  
 $\theta_{new}=\theta_{new}(\theta)$,  
$\phi_{new}=\pm\phi+\phi_0$.  
  
Then, the pull-back of the metric (\ref{metric12}) has the form  
$$(h_1-A_1 \sin\theta_{new}\sin(\pm  
\phi+\phi_0))\left(\left(\frac{\partial \theta_{new}}{\partial  
\theta}\right)^2d\theta^2+\frac{d\phi^2}{4+\cot  
^2(\theta_{new})}\right).
$$ 
If this metric coincides with the  
metric (\ref{metric11}),  
 the coefficients of $d\theta^2$ must  
coincide. Then,  
 $$ \left(\frac{\partial \theta_{new}}{\partial  
\theta}\right)^2(h_1-A_1\sin(\theta_{new})\sin(\pm  
\phi+\phi_0))=h-\frac{A\sin\theta \cos\phi  
}{\sqrt{\cos(\theta)+s}}-\frac{c}{\cos(\theta)+s}, 
$$ 
which implies  (for simplicity, we assume $\textrm{sign}(A)=\textrm{sign}(A_1)$, the case $\textrm{sign}(A)=-\textrm{sign}(A_1)$ can be triated similarly)
\begin{equation}\label{-1}  
A_1\sin(\theta_{new})\left(\frac{\partial \theta_{new}}{\partial  
\theta}\right)^2=\frac{A\sin\theta }{\sqrt{\cos(\theta)+s}}\, ,  
\end{equation}  
\begin{equation}\label{-2} h_1\left(\frac{\partial \theta_{new}}{\partial  
\theta}\right)^2=\left(h-\frac{c}{\cos(\theta)+s}\right). $$  
\end{equation}

Similarly, comparing coefficients of $d\phi^2$, we obtain  
\begin{equation}\label{-3}  
\frac{A_1\sin(\theta_{new})}{4+\cot^2\theta_{new}}=\frac{A\sin(\theta)}{\sqrt{\cos(\theta)+s}}\, ,  
\end{equation}  
\begin{equation}\label{-4}\frac{h_1}{4+\cot^2\theta_{new}}=  
\frac{h-c/(\cos\theta+s)}{1/\sin^2\theta+G(\cos\theta)} \, .  
\end{equation}

Comparing Equations (\ref{-1},\ref{-3}) and Equations  
(\ref{-2},\ref{-4}), we obtain $$\left(\frac{\partial  
\theta_{new}}{\partial \theta}\right)^2= {4+\cot^2\theta_{new}} \  
\ \textrm{and}  \ \ \  \left(\frac{\partial \theta_{new}}{\partial  
\theta}\right)^2=  
\frac{4+\cot^2\theta_{new}}{1/\sin^2\theta+G(\cos\theta)}. $$  
Thus, we obtain a contradiction, since clearly  
${1/\sin^2\theta+G(\cos\theta)}$ is not constant.  
Finally, there exists no diffeomorphism $D$ that takes the metric  
(\ref{metric11}) to metric (\ref{metric12}).

\section{What to do next}  

We constructed a new natural integrable Hamiltonian system on the  
sphere.  It suggests the following ways of further investigation:

  
\begin{enumerate}  
  
\item Construction of bifurcation diagram and description of the  
topological structure of the system.  
  
\item Finding   action-angle variables and explicit formulae for orbits.  
  
\item Quantization of the system.  
  
\item Painlev\'{e} analysis.  

\end{enumerate}

We constructed this new system  using methods based on an  
observation  from  \cite{DMT}. We will describe the methods and  
how we applied them elsewhere.  
It seems that it is the only new integrable system with cubic integral  
on a closed manifold  containing elementary functions only 
that can be found using these methods.
But still one can construct  many new local examples of metrics  
with integrable geodesic flows, which we will do in a joint paper  
with Tabov and Topalov.

\end{document}